\documentclass[10pt]{article}
\usepackage{ngerman,amsmath,amssymb}
\newtheorem{Th}{Theorem}[section]
\newtheorem{Le}{Lemma}[section]
\newtheorem{Co}{Corollary}[section]

\date{}
\begin{document}

\title{On a mixed problem for the nonstationary Stokes system in an angle}
\author{J\"urgen Rossmann}

\maketitle

\section{Introduction}

Let $K$ be an angle with aperture $\alpha$, $K=\{ x=(x_1,x_2): \ 0<r<\infty,\ 0<\varphi<\alpha\}$, where $r$, $\varphi$ denote the polar
coordinates of the point $x$. The sides of $K$ are the half-lines $\Gamma_1$ and $\Gamma_2$,
where $\varphi=0$ on $\Gamma_1$ and $\varphi=\alpha$ on $\Gamma_2$.
The present paper deals with a mixed initial-boundary value problem for the nonstationary Stokes system in $K$,
i. e., we consider the problem
\begin{eqnarray} \label{stokes1a}
&&\frac{\partial {\mathfrak u}}{\partial t} - 2\, \mbox{Div}\, \varepsilon({\mathfrak u}) + \nabla {\mathfrak p} = {\mathfrak f}, \quad
-\nabla\cdot {\mathfrak u} = {\mathfrak g} \ \mbox{ in } Q_K=K\times (0,\infty), \\ \label{stokes2a}
&& {\mathfrak u}(x,t)=0 \ \mbox{ for }x\in \Gamma_1, \ t>0, \  -{\mathfrak p}n+ 2\, \varepsilon_n({\mathfrak u})={\mathfrak h} \
  \mbox{ for }x\in \Gamma_2, \ t>0, \\ \label{stokes3a}
&& {\mathfrak u}(x,0)=0 \ \mbox{ for }x\in K.
\end{eqnarray}
Here, $n$ is the exterior normal to  $\Gamma_2$, $\varepsilon({\mathfrak u})$ denotes the strain tensor with the elements
\[
\varepsilon_{i,j}({\mathfrak u})=\frac 12\,\Big( \frac{\partial{\mathfrak u}_i}{\partial x_j}+\frac{\partial{\mathfrak u}_j}{\partial x_i}\Big), \quad
i,j=1,\ldots,N,
\]
$\mbox{Div}\, \varepsilon({\mathfrak u})$ is the vector with the components $\partial_{x_1}\varepsilon_{i,1} +\cdots +
\partial_{x_N}\varepsilon_{i,N}$, $i=1,\ldots,N$, i.~e., $2\, \mbox{Div}\, \varepsilon({\mathfrak u}) = \Delta {\mathfrak u} +
\nabla\nabla\cdot{\mathfrak u}$, and $\varepsilon_n({\mathfrak u}) = \varepsilon({\mathfrak u})\, n$.

The Dirichlet and the Neumann problems for the Stokes system in an angle were already studied in \cite{kr-16,kr-20,kr-23,r-18}.
As there, the major part of the present paper deals with the parameter-depending problem
\begin{eqnarray} \label{par1}
&& s\, u - \Delta u -\nabla\nabla\cdot u+ \nabla p = f, \quad -\nabla\cdot u = g \mbox{ in } K, \\  \label{par2}
&& u=0 \mbox{ on } \Gamma_1, \quad -pn+ 2\, \varepsilon_n(u)=h    \mbox{ on } \Gamma_2,
\end{eqnarray}
which arises after the Laplace transformation with respect to the time $t$. We show
that this problem has a unique solution $(u,p) \in (V_\beta^2(K)\cap V_\beta^0(K))\times V_\beta^1(K)$ for arbitrary
$f \in V_\beta^0(K)$, $g\in V_\beta^1(K)\cap (V_{-\beta}^1(K,\Gamma_2))^*$, $h\in V_\beta^{1/2}(\Gamma_2)\cap V_\beta^0(\Gamma_2)$
if $\alpha<\pi$ and $\beta$ satisfies the inequalities given in Theorem \ref{bt1}. Furthermore, we obtain the estimate
\begin{eqnarray} \label{1bt1} \nonumber
&& \| u\|_{V_\beta^2(K)} + |s|\, \| u\|_{V_\beta^0(K)} + \| p\|_{V_\beta^1(K)} \le c\, \Big( \| f\|_{V_\beta^0(K)} + \| g\|_{V_\beta^1(K)} \\
&& \quad + |s|\, \| g\|_{V_{-\beta}^1(K,\Gamma_2))^*} + \| h\|_{V_\beta^{1/2}(\Gamma_2)} + |s|^{1/4}\, \| h\|_{V_\beta^0(\Gamma_2)}\Big)
\end{eqnarray}
for this solution, where the constant $c$ is independent of $f,g,h$ and $s$. Here $V_\beta^l(K)$ is the weighted Sobolev space of all functions 
(vector functions) $u$ on $K$ such that $r^{\beta-l+|\alpha|}\partial^\alpha u \in L_2(K)$ for $|\alpha|\le l$.
Using the properties of the Laplace transform, we prove the
existence and uniqueness of a solution $({\mathfrak u},{\mathfrak p})$ of the problem (\ref{stokes1a})--(\ref{stokes3a}), where
${\mathfrak u} \in L_2(0,\infty;V_\beta^2(K))$, $\partial_t{\mathfrak u} \in L_2(0,\infty;V_\beta^0(K))$ and
${\mathfrak p} \in L_2(0,\infty;V_\beta^1(K))$.

\section{The mixed problem for the stationary Stokes system}

We consider the problem (\ref{par1}), (\ref{par2}) in the case $s=0$, i.e.,
\begin{eqnarray} \label{stat1}
&& - \Delta u -\nabla\nabla\cdot u+ \nabla p = f, \quad -\nabla\cdot u = g \mbox{ in } K, \\  \label{stat2}
&& u=0 \mbox{ on } \Gamma_1, \quad -pn+ 2\, \varepsilon_n(u)=h    \mbox{ on } \Gamma_2,
\end{eqnarray}
\underline{Weighted spaces}\\*[1ex]
For nonnegative integer $l$ and real $\beta$, we define the weighted Sobolev spaces $V_\beta^l(K)$
as the sets of all functions (or vector functions) with finite norm
\begin{equation} \label{Vbeta}
\| u\|_{V_\beta^l(K)} = \Big( \int_K r^{2(\beta-l+|\alpha|)}\, \big| \partial_x^\alpha u(x)\big|^2\, dx\Big)^{1/2}
\end{equation}
By $V_\beta^{l-1/2}(\Gamma_j)$, we denote the trace space on $\Gamma_j$ for $V_\beta^l(K)$, $l\ge 1$, which is equipped with the norm
\[
\| u\|_{V_\beta^{l-1/2}(\Gamma_j)} = \inf\Big\{ \| v\|_{V_\beta^l(K)}: \ v\in V_\beta^l(K), \ v=u \mbox{ on }\Gamma_j\Big\}.
\]
Furthermore, $V_\beta^l(K,\Gamma_j)$, $l\ge 1$, is defined as the subspace of all $u\in V_\beta^l(K)$ such that $u=0$ on $\Gamma_j$.\\

\noindent\underline{The operator pencil generated by the mixed problem}\\*[1ex]
Let $\lambda$ be a complex number. We consider the operator
\begin{eqnarray*}
{\cal L}(\lambda):\ \Big( \begin{array}{c} U \\ P \end{array}\Big) \to \left( \begin{array}{c}
  r^{2-\lambda} \big( -\Delta(r^\lambda U(\varphi))-\nabla\nabla\cdot(r^\lambda U(\varphi))+ \nabla (r^{\lambda-1}P(\varphi)\big) \\
  -r^{1-\lambda}\, \nabla\cdot (r^\lambda U(\varphi)) \\
  U(0) \\ -P(\alpha)n +2r^{1-\lambda} \varepsilon_n(r^\lambda U(\varphi)|_\varphi=\alpha \end{array}\right)
\end{eqnarray*}
which maps the space $W^2((0,\alpha))\times W^1((0,\alpha))$ into $L_2((0,\alpha)) \times W^1((0,\alpha))\times {\Bbb C}^2 \times {\Bbb C}^2$.
As is known (see, e.~g., \cite{os-95}) the eigenvalues of the pencil ${\cal L}$ are the solutions of the equation
\begin{equation} \label{ev2}
\lambda^2\, \sin2\alpha - \cos^2(\lambda\alpha)=0.
\end{equation}
Let $\lambda_1(\alpha)$ be the smallest positive solution of the equation (\ref{ev2}). This number is also the solution
with smallest positive real part of the equation (\ref{ev2}). One can easily prove that $\lambda_1(\alpha) > \frac 12$ if $\alpha <\pi$. \\

\noindent\underline{Existence of strong and weak solutions of the problem (\ref{stat1}), (\ref{stat2})}\\*[1ex]
Obviously, the operator $(u,p) \to (f,g,h)$ of the problem (\ref{stat1}), (\ref{stat2}) realizes a linear and continuous mapping
\[
\mbox{from } X_\beta=V_\beta^2(K,\Gamma_1) \times V_\beta^1(K) \quad\mbox{into}\quad
Y_\beta=V_\beta^0(K) \times V_\beta^1(K) \times V_\beta^{1/2}(\Gamma_2).
\]
We denote this operator by $A_\beta$. Besides, one can consider weak (variational) solutions of
the problem (\ref{stat1}), (\ref{stat2}), i.e, pairs $(u,p) \in V_\beta^1(K,\Gamma_1)\times V_\beta^0(K)$ which satisfy
the equations
\begin{equation} \label{varstat}
b_0(u,\overline{v}) - \int p\, \nabla\cdot \overline{v}\, dx = (F,v)_K \ \mbox{for all }v\in V_{-\beta}^1(K,\Gamma_1), \ \
-\nabla\cdot u = g  \mbox{ in }K,
\end{equation}
where
\[
b_0(u,v) = 2 \int_K \sum_{i,j=1}^2 \varepsilon_{i,j}(u)\, \varepsilon_{i,j}(v)\, dx,
\]
$g\in V_\beta^0(K)$, $F$ is a linear and continuous  functional on $V_{-\beta}^1(K,\Gamma_1)$ and $(\cdot,\cdot)_K$ denotes the extension of
the sesquilinear form
\begin{equation}
(F,v)_K = \int_K F \cdot \overline{v}\, dx
\end{equation}
to pairs $(F,v) \in (V_{-\beta}^1(K,\Gamma_1))^* \times V_{-\beta}^1(K,\Gamma_1)$.
The following theorem can be found, e.g., in \cite[Theorem 1.2.5]{mr-10}).

\begin{Th} \label{at1}
\begin{enumerate}
\item $A_\beta$ is an isomorphism if there are no solutions of the equation {\em (\ref{ev2})} on the line $\mbox{\em Re}\, \lambda =\beta-1$.
\item If $(u,p)\in X_\beta$, $A_\beta(u,p) \in Y_\beta\cap Y_\gamma$ and the closed strip between the lines
$\mbox{\em Re}\, \lambda =\beta-1$ and $\mbox{\em Re}\, \lambda =\gamma-1$ is free of solutions of the equation {\em (\ref{ev2})},
then $(u,p) \in X_\gamma$.
\item Suppose that the line $\mbox{\em Re}\, \lambda = \beta$ is free of solutions of the equation {\em (\ref{ev2})}. Then the problem
{\em (\ref{varstat})} has a unique solution $(u,p)\in V_\beta^1(K,\Gamma_1)\times V_\beta^0(K)$ for arbitrary
$F\in (V_{-\beta}^1(K,\Gamma_1))^*$ and $g\in V_\beta^0(K)$. If moreover $g\in V_{\beta+1}^1(K)$ and the functional $F$ has
the form
\[
(F,v)_K = \int_K f\cdot \overline{v}\, dx + \int_{\Gamma_2} h\cdot \overline{v}\, dr \ \mbox{ for all } v\in V_{-\beta}^1(K,\Gamma_1),
\]
where $f\in V_{\beta+1}^0(K)$ and $h\in V_{\beta+1}^{1/2}(\Gamma_2)$, then $(u,p) \in X_{\beta+1}$ and $A_{\beta+1}(u,p)=(f,g,h)$.
\end{enumerate}
\end{Th}

\section{Weak solutions of the parameter-depending problem in the angle \boldmath $K$}

Now we consider the problem (\ref{par1}), (\ref{par2}) with an arbitrary complex parameter $s$.

\subsection{Weighted spaces in an angle}

Let $V_\beta^l(K)$ be the weighted space introduced in Section 2. We set $E_\beta^l(K) = V_\beta^l(K) \cap V_\beta^0(K)$.
Note that the spaces $V_\beta^l(K)$ and $E_\beta^l(K)$ can be also defined as the closures of the set
$C_0^\infty(\overline{K}\backslash \{ 0\})$ with respect to the norms (\ref{Vbeta}) and
\[
\| u\|_{E_\beta^l(K)} =  \| u\|_{V_\beta^l(K)} + \| u\|_{V_\beta^0(K)}\, ,
\]
respectively. An equivalent norm in $E_\beta^l(K)$ is
\[
\| u\| = \Big( \int_K \sum_{|\alpha|\le l} \big( r^{2\beta}+r^{2(\beta-l+|\alpha|)}\big)\big|
  \partial_x^\alpha u(x)\big|^2\, dx\Big)^{1/2}.
\]
For $l\ge 1$, let $E_\beta^l(K,\Gamma_j)=V_\beta^l(K,\Gamma_j) \cap V_\beta^0(K)$, and let $E_\beta^{l-1/2}(\Gamma_j)$ be the trace
space for $E_\beta^l(K)$ an $\Gamma_j$ which is equipped with the norm
\[
\| u\|_{E_\beta^{l-1/2}(\Gamma_j)} = \inf\Big\{ \| v\|_{E_\beta^l(K)}: \ v\in E_\beta^l(K), \ v=u \mbox{ on }\Gamma_j\Big\}.
\]
Note that the norm
\begin{equation} \label{tracenorm}
\| u\| = \| u\|_{V_\beta^{1/2}(\Gamma_j)} + \| u\|_{V_\beta^0(\Gamma_j)}\, ,\ \mbox{ where }
  \| u\|_{V_\beta^0(\Gamma_j)} = \| r^\beta u\|_{L_2(\Gamma_j)},
\end{equation}
is equivalent to the $E_\beta^{1/2}(\Gamma_j)$-norm and that $V_\beta^{1/2}(\Gamma_j) \subset V_{\beta-1/2}^0(\Gamma_j)$,
where $V_\beta^0(\Gamma_j)$ is the set of all functions $u$ on $\Gamma_j$ such that $r^\beta u \in L_2(\Gamma_j)$
(cf. \cite[Lemmas 1.4 and 1.6]{mp-78} or \cite[Lemmas 6.1.2 and 6.5.4]{kmr-97}). Hence,
$E_\beta^{1/2}(\Gamma_j)\subset V_\beta^0(\Gamma_j)\cap V_{\beta-1/2}^0(\Gamma_j)$.

The dual spaces of $V_\beta^l(K)$ and $V_\beta^l(K,\Gamma_j)$ are denoted by $(V_\beta^l(K))^*$ and
$(V_\beta^l(K,\Gamma_j))^*$, respectively. Obviously, the sesquilinear form
\begin{equation} \label{sprod}
(u,v)_K = \int_K u\cdot \bar{v}\, dx
\end{equation}
is continuous on $V_\beta^0(K)\times V_{-\beta}^0(K)$. Hence, $(V_\beta^0(K))^* = V_{-\beta}^0(K)$ and
$V_{l-\beta}^0(K) \subset (V_\beta^l(K))^* \subset (V_\beta^l(K,\Gamma_j))^*$ for $l\ge 1$.

\subsection{Solutions of the equation \boldmath $\mbox{div}\, u = g$}

Obviously, the operator div continuously maps $V_\beta^1(K)$ into $V_\beta^0(K)$. If $u\in V_\beta^1(K,\Gamma_1)\cap V_\beta^0(K)$ and
$v\in V_{-\beta}^1(K,\Gamma_2)$, then
\[
\Big| \int_K (\nabla\cdot u)\, v\,dx\Big| = \Big| \int_K u\cdot \nabla v\, dx\Big| \le \| u\|_{V_\beta^0(K)}\ \| v\|_{V_{-\beta}^1(K)}
\]
Thus, the operator div realizes a linear and continuous mapping from $E_\beta^1(K,\Gamma_1)$ into
$V_\beta^0(K) \cap (V_{-\beta}^1(K,\Gamma_2))^*$. We prove that this operator is surjective.

\begin{Le} \label{al1}
Suppose that $g \in V_\beta^0(K) \cap (V_{-\beta}^1(K,\Gamma_2))^*$.
Then there exists a vector function $u\in E_\beta^1(K,\Gamma_1)$ satisfying the equation
$-\nabla\cdot u  = g$ and the estimate
\begin{equation} \label{1al1}
\| u\|_{E_\beta^1(K)} \le c\, \Big( \| g\|_{V_\beta^0(K)} + \| g\|_{(V_{-\beta}^1(K,\Gamma_2))^*}\Big)
\end{equation}
with a constant $c$ independent of $g$.
\end{Le}

Proof.
1. First, we assume that the angle $\alpha$ is sufficiently small such that
that $|\beta|$ and $|\beta-1|$ are less than $\pi/\alpha$ and the strip $\beta \le \mbox{Re}\, \lambda \le \beta+1$ does not contain
solutions of the equation (\ref{ev2}). Since $g\in V_\beta^0(K) \cap (\stackrel{\circ}{V}\!{}_{-\beta}^1(K))^*$ and the
interval $[\beta-1,\beta]$ does not contain eigenvalues $j\pi/\alpha$ of the operator pencil generated by the Dirichlet problem for the Laplacian,
there exists a unique solution $u_0 \in V_\beta^2(G)\cap \stackrel{\circ}{V}\!{}_\beta^1(K)$ of the equation $\Delta u_0=g$
in $K$. Then $U=-\nabla u_0$ satisfies the equation $-\nabla\cdot U=g$ and the estimates
\[
\| U\|_{V_\beta^1(K)} \le c\, \| g\|_{V_\beta^0(K)}\, , \quad \| U\|_{V_\beta^0(K)}  \le c\, \| g\|_{(\stackrel{\circ}{V}{}_{-\beta}^1(K))^*}
  \le c \, \| g\|_{(V_{-\beta}^1(K,\Gamma_2))^*}\, .
\]
However, the vector function $U$ does not necessarily satisfy the condition $U=0$ on $\Gamma_1$. Let $\chi_1$ be a continuously differentiable function
on the interval $[0,\alpha]$ such that $\chi_1(0)=1$ and $\chi_1(\alpha)=0$.
If $r,\varphi$ are the polar coordinates of the point $x$, we can consider the function $\chi_1$, i.e. the mapping $x\to \chi_1(\varphi)$,
as a function on $\overline{K}\backslash \{ 0\}$. We define $\chi_2=1-\chi_1$. Then $\chi_2 U \in E_\beta^1(K,\Gamma_1)$ and
\[
- \nabla\cdot (\chi_2 U)  = \chi_2 g + U\cdot \nabla\chi_1.
\]
We construct a vector function $w\in V_\beta^1(K,\Gamma_1)\cap V_\beta^0(K)$ such that $-\nabla\cdot w =\chi_1 g - U\cdot\nabla\chi_1$.
Since $\chi_1=0$ on $\Gamma_2$ and $|\nabla\chi_1| \le c\, r^{-1}$, it follows that $\chi_1 g\in (V_{-\beta}^1(K))^*$ and
\[
\| \chi_1 g\|_{(V_{-\beta}^1(K))^*} \le c\, \| g\|_{(V_{-\beta}^1(K,\Gamma_2))^*}\, .
\]
Furthermore,
\[
\| U\cdot\nabla\chi_1 \|_{(V_{-\beta}^1(K))^*} \le c\, \| U\|_{V_\beta^0(K)} \ \mbox{ and } \
\| U\cdot\nabla\chi_1 \|_{V_\beta^0(K)} \le c\, \| U\|_{V_\beta^1(K)}\, .
\]
By Theorem \ref{at1}, there exists a solution of $(w,p)\in V_\beta^1(K,\Gamma_1)\times V_\beta^0(K)$ of the problem
\[
b_0(w,v) - \int_K p\, \nabla\cdot v\, dx = 0 \ \mbox{ for all }v\in V_{-\beta}^1(K,\Gamma_1), \quad -\nabla\cdot w
  =\chi_1 g - U\cdot\nabla\chi_1 \ \mbox{ in }K
\]
satisfying the estimate
\[
\| w\|_{V_\beta^1(K)} + \| p\|_{V_\beta^0(K)} \le c\, \| \chi_1 g - U\cdot \nabla\chi_1\|_{V_\beta^0(K)}\, .
\]
This implies
\[
\int_K w\cdot (-\Delta v-\nabla\nabla\cdot v+ \nabla q)\, dx + \int_K (q\, \nabla\cdot w- p\nabla\cdot v)\, dx
+ \int_{\Gamma_2} u\cdot (-qn+2\varepsilon_n(v))\, dr =0
\]
for all $v\in V_{1-\beta}^2(K,\Gamma_1)$ and $q\in V_{1-\beta}^1(K)$. Let $\phi$ be a vector function in
$V_{-\beta}^0(K)\cap V_{1-\beta}^0(K)$. By Theorem \ref{at1}, there exists a unique solution
$(v,q)$, $v\in V_{-\beta}^2(K)\cap V_{1-\beta}^2(K)$, $q\in V_{-\beta}^1(K)\cap V_{1-\beta}^1(K)$, of the system
\[
-\Delta v-\nabla\nabla\cdot v+ \nabla q = \phi, \ \nabla\cdot v=0\ \mbox{ in }K,
\]
with the boundary conditions $v=0$ on $\Gamma_1$, $-qn+2\varepsilon_n(v)=0$ on $\Gamma_2$. Furthermore,
\[
\| v\|_{V_{-\beta}^2(K)} + \| q\|_{V_{-\beta}^1(K)} \le c\, \| \phi\|_{V_{-\beta}^0(K)}\, .
\]
Then we obtain
\begin{eqnarray*}
\Big| \int_K w\cdot \phi\, dx\Big| & = & \Big| \int_K (\chi_1 g - U\cdot \nabla\chi_1)\, q\, dx\Big|
\le \| \chi_1 g - U\cdot \nabla\chi_1\|_{(V_{-\beta}^1(K))^*}\, \| q\|_{V_{-\beta}^1(K))} \\
& \le & c\, \| \chi_1 g - U\cdot \nabla\chi_1\|_{(V_{-\beta}^1(K))^*}\, \| \phi\|_{V_{-\beta}^0(K))}
\end{eqnarray*}
for all $\phi \in V_{-\beta}^0(K)\cap V_{1-\beta}^0(K)$. Consequently,
\[
\| w\|_{V_\beta^0(K)} \le c\,  \| \chi_1 g - U\cdot \nabla\chi_1\|_{(V_{-\beta}^1(K))^*} \le c'\, \| g\|_{(V_{-\beta}^1(K,\Gamma_2))^*}\, .
\]
Thus $u=\chi_2 U + w$ satisfies the equation $-\nabla\cdot u =g$, the condition $u=0$ on $\Gamma_1$ and the estimate (\ref{1al1}).
This proves the lemma for small $\alpha$.

2. If $\alpha$ is arbitrary less than $2\pi$, then we can argue as in the proof of \cite[Lemma 2.6]{r-18} and \cite[Lemma 3.3]{kr-23}.
The equation $\nabla\cdot u = -g$ can be written as
\begin{equation} \label{2ac2}
\partial_r\, u_r + r^{-1}\, \big( u_r + \partial_\varphi u_\varphi\big)=-g(r,\varphi) \quad\mbox{for }r>0,\ 0<\varphi<\alpha.
\end{equation}
with the polar components $u_r,u_\varphi$. Let $\delta$ be a small
positive number and let $v_r(r,\varphi)=u_r(r,\delta^{-1}\varphi)$, $v_\varphi(r,\varphi)=\delta\, u_\varphi(r,\delta^{-1}\varphi)$ and
$h(r,\varphi)= g(r,\delta^{-1}\varphi)$ for $r>0$ and $0<\varphi<\delta\alpha$. Then (\ref{2ac2}) is equivalent to
\begin{equation} \label{3ac2}
\partial_r\, v_r + r^{-1}\, \big( v_r + \partial_\varphi v_\varphi\big)=-h(r,\varphi) \quad\mbox{for }r>0,\ 0<\varphi<\delta\alpha.
\end{equation}
By the first part of the proof, the equation (\ref{3ac2}) has a solution $(v_r,v_\varphi)$ satisfying the estimate
\begin{eqnarray*}
&&\int_0^\infty\int_0^{\delta\alpha} r^{2\beta-1} \Big( (1+r^2)\, \big| (v_r,v_\varphi)\big|^2
  + \big|\partial_\varphi(v_r,v_\varphi)\big|^2 + r^2\, \big|\partial_r (v_r,v_\varphi)\big|^2\Big)\, d\varphi\, dr \\
&& \le c\,  \Big( \| h\|^2_{V_\beta^0(K')} + \| h\|^2_{V_{-\beta}^1(K',\Gamma'_N))^*}\Big),
\end{eqnarray*}
where $K'$ is the angle $0<\varphi<\delta\alpha$ and $\Gamma'_2$ is the half-line $\varphi=\delta\alpha$.
Then $u_r$ and $u_\varphi$ satisfy (\ref{2ac2}). Furthermore, the estimate (\ref{1al1}) holds. $\Box$

\subsection{Existence of weak solutions}

Let $b_s$ denote the bilinear form
\[
b_s(u,v) = b_0(u,v)+ \int_K s u\cdot v\, dx = \int_K \Big( su\cdot v + 2\sum_{i,j=1}^2 \varepsilon_{i,j}(u)\, \varepsilon_{i,j}(v)\Big)\, dx.
\]
Obviously,
\begin{equation} \label{estbs}
|b_s(u,v)|\le c\, \max(1,|s|)\,  \| u\|_{E_\beta^1(K)} \, \| v\|_{E_{-\beta}^1(K)}
\end{equation}
for arbitrary $u\in E_\beta^1(K)$, $v\in E_{-\beta}^1(K)$, where $c$ is independent of $u,v$ and $s$.
Analogously to \cite[Lemma 3.1]{kr-23}, the following lemma holds.

\begin{Le} \label{al2}
Suppose that $(u,p) \in E_\gamma^2(K,\Gamma_1) \times V_\gamma^1(K)$ and $p|_{\Gamma_2}\in E_\gamma^{1/2}(\Gamma_2)$, where
$\beta \le \gamma \le \beta+1$. Then the following assertions are true.

1) $u\in E_\beta^1(K,\Gamma_1)$, $p\in V_\beta^0(K)+V_\beta^1(K,\Gamma_2)$ and
\[
\| u\|_{E_\beta^1(K)} + \| p\|_{V_\beta^0(K)+V_\beta^1(K,\Gamma_2)}
  \le c\, \Big( \| u\|_{E_\gamma^2(K)} + \| p\|_{V_\gamma^1(K)} + \| p\|_{E_\gamma^{1/2}(\Gamma_2)}\Big)
\]
with a constant $c$ independent of $u$ and $p$.

2) $u$ and $p$ satisfy the equality
\begin{equation} \label{1al2}
b_s(u,\bar{v}) - \int_K  p\, \nabla \cdot \bar{v}\, dx  = (F,v)_K \quad\mbox{for all } v\in E_{-\beta}^1(K,\Gamma_1),
\end{equation}
where
\[
(F,v)_K = \int_K \big( su-\Delta u -\nabla\nabla\cdot u + \nabla p\big)\cdot \bar{v}\, dx
 + \int_{\Gamma_2} \big( -pn+2\varepsilon_n(u)\big)\cdot \bar{v}\, dr
\]
and
\begin{equation} \label{2al2}
\| F \|_{(E_{-\beta}^1(K,\Gamma_1))^*} \le c\, \Big( \| u\|_{E_\gamma^2(K)} + \| p\|_{V_\gamma^1(K)} + \| p\|_{E_\gamma^{1/2}(\Gamma_2)} \Big)
\end{equation}
with a constant $c$ independent of $u$ and $p$.
\end{Le}

Proof.
1) Obviously, $E_\gamma^2(K) \subset E_\beta^1(K)$ for $\beta\le \gamma\le \beta+1$. Let $\zeta$ be a smooth (two times continuously differentiable) function on ${\mathbb R}_+$ with support in $[0,1]$ which is equal to one on the interval $(0,\frac 12)$, and let $\eta=1-\zeta$. The functions $\zeta,\eta$ can be considered as smooth functions on $\overline{K}$ if one defines $\zeta(x)=\zeta(|x|)$ and $\eta(x)=\eta(|x|)$.
If $p\in V_\gamma^1(K)$ and $p|_{\Gamma_2}\in E_\beta^{1/2}(\Gamma_2)$, then there exists a function
$q\in E_\gamma^1(K)$ such that $q=p$ on $\Gamma_2$, i.e., $p-q\in V_\gamma^1(K,\Gamma_2)$.
Since $\gamma-1 \le \beta \le \gamma$, it follows that $\zeta(p-q) \in V_\beta^0(K)$ and $\eta(p-q) \in V_\beta^1(K,\Gamma_2)$.
Furthermore $q \in E_\gamma^1(K)\subset V_\beta^0(K)$. This means $p=p_1+p_2$,
where $p_1=\zeta p+\eta q = q+\zeta(p-q)\in V_\beta^0(K)$ and $p_2=\eta(p-q)\in V_\beta^1(K,\Gamma_2)$. Moreover,
\[
\| p\|_{V_\beta^0(K)+V_\beta^1(K,\Gamma_2)} \le \| p_1\|_{V_\beta^0(K)} + \| p_2 \|_{V_\beta^1(K)}
  \le c\, \Big( \| p\|_{V_\gamma^1(K)} + \| p\|_{E_\gamma^{1/2}(\Gamma_2)}\Big) .
\]
2) We estimate the norm of the functional $F$ in $E_{-\beta}^1(K,\Gamma_2))^*$. Since $E_{-\beta}^1(K)\subset V_{-\gamma}^0(K)$
for $\beta\le \gamma\le \beta+1$, we have
\[
su - \Delta u - \nabla\nabla\cdot u + \nabla p \in V_\gamma^0(K) \subset (E_{-\beta}^1(K))^* \subset (E_{-\beta}^1(K,\Gamma_2))^*.
\]
Furthermore, $h= -pn + 2\varepsilon_n(u) \in E_\gamma^{1/2}(\Gamma_2)$ if $u\in E_\gamma^2(K)$ and $p|_{\Gamma_2}\in E_\gamma^{1/2}(\Gamma_2)$.
Consequently,
\begin{eqnarray*}
&& \Big| \int_{\Gamma_2} h\cdot v\, d\sigma\Big|^2 \le \int_{\Gamma_2} (r^{2\gamma-1} + r^{2\gamma})\, |h|^2\, d\sigma \
\int_{\Gamma_2} (r^{-2\beta-1} + r^{-2\beta})\, |v|^2\, d\sigma \\
&& \le \| h\|^2_{E_\gamma^{1/2}(\Gamma_2)} \ \| v\|^2_{E_{-\beta}^{1/2}(\Gamma_2)} \le
  \| h\|^2_{E_\gamma^{1/2}(\Gamma_2)} \ \| v\|^2_{E_{-\beta}^1(K)}
\end{eqnarray*}
for $v\in E_{-\beta}^1(K)$, $\beta\le \gamma\le \beta+1$. This implies (\ref{2al2}).
Integrating by parts, the formula (\ref{1al2}) holds for $v\in C_0^\infty(\overline{K}\backslash\{ 0\})$.
The right-hand side of (\ref{1al2}) satisfies the estimate (\ref{2al2}). Furthermore, the inequality (\ref{estbs})
is valid for $u\in E_\beta^1(K)$, while the estimate
\[
\Big| \int_K p\, \nabla\cdot v\, dx\Big| = \Big| \int_K \big( p_1\nabla\cdot v - v\cdot \nabla p_2 \big)\, dx
  \le \Big( \| p_1\|_{V_\beta^0(K)} + \| p_2\|_{V_\beta^1(K)} \Big)\, \| v\|_{E_{-\beta}^1(K)}
\]
holds for arbitrary $p=p_1+p_2 \in V_\beta^0(K) + V_\beta^1(K,\Gamma_2)$ and $v\in E_{-\beta}^1(K,\Gamma_1)$.
Since the set $C_0^\infty(\overline{K}\backslash\{ 0\})$ is dense in $E_{-\beta}^1(K)$, it follows that (\ref{1al1}) is valid for
arbitrary $v\in E_{-\beta}^1(K)$. The proof is complete. $\Box$ \\

Motivated by Lemma \ref{al2}, we define weak solutions of the problem (\ref{par1}), (\ref{par2}) as follows.
Let $F\in (E_{-\beta}^1(K,\Gamma_1))^*$ and $g\in V_\beta^0(K)\cap (V_{-\beta}^1(K,\Gamma_2))^*$ be given. The vector function
$(u,p)\in E_\beta^1(K,\Gamma_1)\times\big( V_\beta^0(K)+V_\beta^1(K,\Gamma_2)\big)$ is called a weak solution
if $(u,p)$ satisfies the identity (\ref{1al2}) for all $v\in E_{-\beta}^1(K,\Gamma_1)$ and the equation
\begin{equation} \label{3al2}
-\nabla\cdot u=g \ \mbox{ in }  K.
\end{equation}
Using Lemma \ref{al1} and the inequality
\begin{equation} \label{bs}
\big| b_s(u,\bar{u})\big| \ge c\, \| u\|^2_{E_0^1(K)} \ \mbox{ for all } u\in E_0^1(K,\Gamma_1), \ \mbox{Re}\, s \ge 0, \ s\not=0
\end{equation}
with a constant $c$ depending on $|s|$ (see \cite[Lemma3.5]{kr-23}), we can prove the existence and uniqueness of weak solutions in the space
$E_0^1(K,\Gamma_1)\times\big( L_2(K)+V_0^1(K,\Gamma_2)\big)$ analogously to \cite[Theorem 1.1]{kr-16} and \cite[Theorem 3.6]{kr-23}.

\begin{Th} \label{at2}
Suppose that $F\in (E_0^1(K,\Gamma_1))^*$, $g\in L_2(K)\cap (V_0^1(K,\Gamma_2))^*$, $\mbox{\em Re}\, s \ge 0$ and $s\not=0$.
Then there exists a uniquely determined solution
$(u,p) \in E_0^1(K,\Gamma_1) \times \big( L_2(K)+ V_0^1(K,\Gamma_2)\big)$
of the problem {\em (\ref{1al2}), (\ref{3al2})} with $\beta=0$. Furthermore,
\[
\| u\|_{E_0^1(K)} + \| p\|_{L_2(K)+V_0^1(K,\Gamma_2)} \le c\, \Big( \| F\|_{(E_0^1(K,\Gamma_1))^*}
+ \| g\|_{L_2(K)} + \| g\|_{(V_0^1(K,\Gamma_2))^*}\Big)
\]
with a constant $c$ depending on $|s|$ but not on $F$ and $g$.
\end{Th}

Proof.
Lemma \ref{al1} allows us to restrict ourselves in the proof to the case $g=0$.
Let ${\cal H} =\{ u\in E_0^1(K,\Gamma_1): \ \nabla\cdot u =0 \ \mbox{in }K\}$
and ${\cal H}^\bot$ its orthogonal complement in $E_0^1(K,\Gamma_1)$. Since the inequality (\ref{bs}) is satisfied for all
$u\in {\cal H}$, there exists a uniquely determined vector function $u\in {\cal H}$ such that
\begin{equation} \label{2at2}
b_s(u,\bar{v})= ( F,v)_K \ \mbox{ for }v\in {\cal H}, \qquad
  \| u\|_{E_0^1(K)} \le c\, \| F\|_{{\cal H}^*} \le c\, \| F\|_{(E_0^1(K,\Gamma_1))^*}\, .
\end{equation}
By Lemma \ref{al1}, the operator $\mbox{div}$ is an isomorphism from ${\cal H}^\bot$ onto $L_2(K)\cap (V_0^1(K,\Gamma_2))^*$,
i.e., for every $q \in L_2(K)\cap (V_0^1(K,\Gamma_2))^*$, there exists a uniquely determined vector function
$v=(-\mbox{div})^{-1}q \in {\cal H}^\bot$ such that $\nabla\cdot v=-q$ and
\[
\| v\|_{E_0^1(K)} \le c\, \Big( \| q\|_{L_2(K)} + \| q\|_{(V_0^1(K,\Gamma_2))^*}\Big) .
\]
We consider the functional
\begin{equation} \label{3at2}
\ell(q) =  (F,\bar{v})_K - b_s(u,v) = (F, (-\mbox{div})^{-1}\bar{q})_K - b_s\big(u,(-\mbox{div})^{-1}q\big)
\end{equation}
on $L_2(K)\cap (V_0^1(K,\Gamma_2))^*$. Obviously,
\[
\big| \ell(q)\big|  \le  c\, \| F\|_{(E_0^1(K,\Gamma_1))^*} \, \| v\|_{E_0^1(K)} \le  c\, \| F\|_{(E_0^1(K,\Gamma_1))^*}\,
  \Big( \| q\|_{L_2(K)} + \| q\|_{(V_0^1(K,\Gamma_2))^*}\Big).
\]
This means that $\ell$ is continuous on $L_2(K)\cap (V_0^1(K,\Gamma_2))^*$ and there exists an element $p$ of the dual
space $L_2(K)+V_0^1(K,\Gamma_2)$ such that $\ell(q) =  (p,\bar{q})_K$ for all $q\in L_2(K)\cap (V_0^1(K,\Gamma_2))^*$.
Then by (\ref{3at2}), we have
\[
b_s(u,\bar{v}) + (p,q)_K = (F,v)_K
\]
for $q\in L_2(K)\cap V_0^{-1}(K)$, $v=(-\mbox{div})^{-1}q \in {\cal H}^\bot$. This means that
\[
b_s(u,\bar{v}) -  \int_K p \, \nabla \cdot \bar{v}\, dx =  (F,v)_K \ \mbox{ for all }v\in {\cal H}^\bot.
\]
Since $b_s(u,\bar{v})=(F,v)_K$ for all $v\in {\cal H}$ (see (\ref{2at2})), it follows that $(u,p)$ is a solution of the problem
(\ref{1al2}), (\ref{3al2}) for $g=0$ and $\beta=0$. The uniqueness of the solution can be easily deduced from (\ref{bs}) and
Lemma \ref{al1}. $\Box$

\begin{Co} \label{ac3}
Suppose that $s\not=0$, $\mbox{\em Re}\, s \ge 0$, $f\in V_\beta^0(K)$, $g\in V_\beta^1(K)\cap (V_{-\beta}^1(K,\Gamma_2))^*$ and
$h\in E_\beta^{1/2}(\Gamma_2)$, where $0\le \beta\le 1$.
Then there exists a unique solution $(u,p) \in E_0^1(K,\Gamma_1) \times \big( L_2(K)+ V_0^1(K,\Gamma_2)\big)$
of the problem {\em (\ref{1al2}), (\ref{3al2})}, where
\begin{equation} \label{1ac1}
(F,v)_K = \int_K f\cdot \bar{v}\, dx + \int_{\Gamma_2} h\cdot \bar{v}\, dr \ \mbox{for all } v\in E_0^1(K,\Gamma_1).
\end{equation}
\end{Co}

Proof. Obviously, $g\in V_\beta^1(K)\cap (\stackrel{\circ}{V}\!{}_{-\beta}^1(K))^*$, where 
$\stackrel{\circ}{V}\!{}_{-\beta}^1(K) = V_{-\beta}^1(K,\Gamma_1) \cap V_{-\beta}^1(K,\Gamma_2)$. 
By  \cite[Lemma 3.3]{kr-23}, there exists a vector function $v\in E_\beta^2(K)$ such that $\nabla\cdot v =g$. 
Thus, $g \in E_\beta^1(K) \subset L_2(K)$ since $0\le \beta\le 1$. Obviously, $V_\beta^1(K) \subset V_{\beta-1}^0(K) \subset (V_{2-\beta}^1(K,\Gamma_2))^*$
and, consequently, $g \in (V_{2-\beta}^1(K,\Gamma_2))^* \cap (V_{-\beta}^1(K,\Gamma_2))^* \subset (V_0^1(K,\Gamma_2))^*$.
We show that the functional (\ref{1ac1}) is continuous on $E_0^1(K,\Gamma_1)$.
Since $E_0^1(K) \subset V_{-\beta}^0(K)$ for $0\le \beta \le 1$
it follows that $f\in V_\beta^0(K) \subset (E_0^1(K,\Gamma_1))^*$. As a consequence of the imbedding $E_\beta^{1/2}(\Gamma_2)\subset
V_\beta^0(\Gamma_2) \cap V_{\beta-1/2}^0(\Gamma_2)$ (see Subsection 3.1), we obtain the estimate
\[
\Big|  \int_{\Gamma_2} h\cdot v\, dr\Big|  \le  \| r^\gamma h\|_{L_2(\Gamma_2)} \ \| r^{-\gamma} v\|_{L_2(\Gamma_2)}
 \le  c\, \| h\|_{E_\beta^{1/2}(\Gamma_2)} \ \| v\|_{E_0^{1/2}(\Gamma_2)}
\]
if $\max(0,\beta-\frac 12) \le \gamma \le \min(\beta,\frac 12)$. Consequently, $F$ is continuous on $E_0^1(K,\Gamma_1)$,
and the assertion of the corollary follows immediately from Theorem \ref{at2}. $\Box$

\section{Strong solutions of the parameter-dependent problem}

Suppose that $s\not=0$ and $u\in E_\beta^2(K,\Gamma_1)$. Then
\[
su-\Delta u -\nabla\nabla\cdot u\in V_\beta^0(K), \ \  \nabla\cdot u \in V_\beta^1(K) \cap (V_{-\beta}^1(K,\Gamma_2))^*
\]
and $\varepsilon_n(u) \in E_\beta^{1/2}(\Gamma_2)$. In order to have $\nabla p \in V_\beta^0(K)$
and $h=-pn+2\varepsilon_n(u) \in E_\beta^{1/2}(\Gamma_2)$, we assume that $p\in V_\beta^1(K)$ and $p|_{\Gamma_2} \in
E_\beta^{1/2}(\Gamma_2)$. The space of such functions $p$ is denoted by ${\cal V}_\beta^1(K)$.
Since $E_\beta^{1/2}(\Gamma_2)= V_\beta^{1/2}(\Gamma_2) \cap V_\beta^0(\Gamma_2)$ (see Subsection 2.1), we provide
the space ${\cal V}_\beta^1(K)$ with the norm
\[
\| p\|_{{\cal V}_\beta^1(K)} = \| p\|_{V_\beta^1(K)} + \| p|_{\Gamma_2} \|_{V_\beta^0(\Gamma_2)}\, .
\]
We denote the operator
\begin{eqnarray} \label{operator}
&& E_\beta^2(K,\Gamma_1) \times {\cal V}_\beta^1(K) \ni (u,p) \nonumber \\
&& \quad \to \ (f,g,h) \in V_\beta^0(K) \times \big( V_\beta^1(K)\cap (V_{-\beta}^1(K,\Gamma_2))^*\big)\times E_\beta^{1/2}(\Gamma_2)
\end{eqnarray}
of the problem (\ref{par1}), (\ref{par2}) by $A_{\beta,s}$.

\subsection{An estimate for the function \boldmath $p$}

Let $(u,p)\in E_\beta^2(K,\Gamma_1)\times {\cal V}_\beta^1(K)$ be a solution of the problem (\ref{par1}), (\ref{par2}), where $s\not=0$.
Then
\begin{equation} \label{problemp}
\int_K \nabla p\cdot \nabla \bar{q}\, dx= (\Phi,q)_K \mbox{ for all } q\in V_{-\beta}^1(K,\Gamma_2),
\quad p = H \ \mbox{ on }\Gamma_2,
\end{equation}
where $\Phi\in (V_{-\beta}^1(K,\Gamma_2))^*$ is defined as
\[
(\Phi,q)_K = \int_K \big( (f-\nabla g+\Delta u)\cdot \nabla \bar{q} -s\, g\bar{q}\big)\, dx \quad\mbox{for all }q\in V_{-\beta}^1(K,\Gamma_2)
\]
and $H=\big( 2\varepsilon_n(u) -h)\cdot n \in E_\beta^{1/2}(\Gamma_2)$.

We assume that $\beta\not=(k+\frac 12)\frac\pi\alpha$ for integer $k$ (i.e., $\beta$ is not an eigenvalue of the operator pencil generated by the
mixed problem for the Laplacian in $K$). Then the  problem (\ref{problemp}) has a unique solution $p \in V_\beta^1(K)$ for arbitrary
$\Phi \in (V_{-\beta}^1(K,\Gamma_2))^*$ and $H \in V_\beta^{1/2}(\partial K)$ (see, e. g., \cite[Theorem 3.3.9]{mr-10}).
Clearly, this solution belongs to the subspace ${\cal V}_\beta^1(K)$ of $V_\beta^1(K)$ if
$H \in E_\beta^{1/2}(\Gamma_2)$.

Obviously, $\Phi=\Phi_1+\Phi_2$, where $\Phi_1,\Phi_2$ denote the functionals
\[
(\Phi_1,q)_K = \int_K \big( (f-2\nabla g)\cdot \nabla \bar{q} -sg\bar{q}\big)\, dx, \quad
(\Phi_2,q)_K = \int_K (\Delta u - \nabla\nabla\cdot u)\cdot \nabla \bar{q}\, dx
\]
on  $V_{-\beta}^1(K,\Gamma_2)$. The function $p$ admits the decomposition $p= p_1 + p_2$,
where $p_1$, $p_2$ are the unique solutions of the problems
\begin{equation} \label{problemp1}
\int_K \nabla p_1 \cdot \nabla \bar{q}\, dx = (\Phi_1,q)_K \ \mbox{ for all } q\in V_{-\beta}^1(K,\Gamma_2),\quad
p_1= -h\cdot n \mbox{ on }\Gamma_2
\end{equation}
and
\begin{equation} \label{problemp2}
\int_K \nabla p_2 \cdot \nabla \bar{q}\, dx = (\Phi_2,q)_K \ \mbox{ for all } q\in V_{-\beta}^1(K,\Gamma_2),\quad
p_2= 2\varepsilon_n(u)\cdot n \mbox{ on }\Gamma_2\, ,
\end{equation}
respectively. If  $\beta\not=(k+\frac 12)\frac\pi\alpha$ for integer $k$, then $p_1$ and $p_2$ satisfy the estimates
\begin{equation} \label{estp1}
\| p_1\|_{V_\beta^1(K)} \le  c\, \Big( \| f\|_{V_\beta^0(K)} + \| g\|_{V_\beta^1(K)} + |s|\, \| g\|_{V_\beta^{-1}(K)}
  + \| h\|_{V_\beta^{1/2}(\partial K)}\Big)
\end{equation}
and
\begin{equation} \label{estp2}
\| p_2\|_{V_\beta^1(K)} \le c'\, \Big( \| \Phi_2\|_{V_\beta^{-1}(K)} + \| \varepsilon_n(u)\|_{V_\beta^{1/2}(\partial K)}\Big)
  \le c\, \| Du\|_{V_\beta^1(K)}
\end{equation}
with a constant $c$ independent of $u,p$ and $s$.

\begin{Le} \label{al3}
Suppose that $s\not=0$ and that $(u,p) \in E_\beta^2(K,\Gamma_1) \times {\cal V}_\beta^1(K)$ is a solution of the problem
{\em (\ref{par1}), (\ref{par2})}. Assume furthermore that $\beta\le \gamma\le \beta+\frac 12$ and that
all numbers $(k+\frac 12)\frac\pi\alpha$ with integer $k$ lie outside the interval $\beta \le \lambda \le \gamma$. Then
$p=p_1+p_2$, where $p_1 \in V_\beta^1(K)$, $p_2 \in V_\beta^1(K)\cap V_{\gamma-1}^0(K)$, $p_1$ and $p_2$ satisfy the estimates
{\em (\ref{estp1})} and {\em (\ref{estp2})}. In addition,
\begin{equation} \label{1al3}
\| p_2\|_{V_{\gamma-1}^0(K)} \le c\,  \| u\|_{E_\beta^2(K)}
  \end{equation}
with a constant $c$ independent of $u$ and $p$. If $f\in V_\gamma^0(K)$, $g\in V_\gamma^1(K)\cap
(V_{-\gamma}^1(K,\Gamma_2))^*$ and $h\in V_\gamma^{1/2}(\Gamma_2)$, then $p\in V_{\gamma-1}^0(K)$.
\end{Le}

Proof. Let $p_1\in V_\beta^1(K)$, $p_2\in V_\beta^1(K)$ be the solutions of the problems (\ref{problemp1}) 
and (\ref{problemp2}), respectively. Then
\[
- \int_K p_2\, \Delta\bar{q}\, dx + \int_\Gamma p_2\, \frac{\partial\bar{q}}{\partial n}\, dr =
  (\Phi_2,q)_\Omega = \int_\Gamma \sum_{i,j=1}^2 \frac{\partial u_i}{\partial x_j} \,
  \Big(  n_j\, \frac{\partial \bar{q}}{\partial x_i} - n_i\, \frac{\partial \bar{q}}{\partial x_j}\Big)\, dr
\]
for all $q\in V_{1-\beta}^2(K,\Gamma_2)$, where $\Gamma=\Gamma_1\cup \Gamma_2$. Let $\phi \in V_{1-\beta}^0(K)\cap V_{1-\gamma}^0(K)$.
Since the interval $\beta \le \lambda \le \gamma$ contains no eigenvalues of the pencil generated by the mixed problem for the
Laplace equation, there exists a solution $q \in V_{1-\beta}^2(K)\cap V_{1-\gamma}^2(K)$ of the problem
$-\Delta q = \phi$ in $K$, $q=0$ on $\Gamma_2$, $\frac{\partial q}{\partial n}=0$ on $\Gamma_1$ satisfying the estimate
\[
\| q\|_{V_{1-\gamma}^2(K)} \le c\, \| \phi\|_{V_{1-\gamma}^0(K)}\, .
\]
Then we obtain
\[
\int_K p_2\, \bar{\phi}\, dx = (\Phi_2,q)_K - \int_{\Gamma_2}\, 2\varepsilon_u(u)\, \cdot n\, \frac{\partial\bar{q}}{\partial n}\, dr.
\]
One can easily show (see \cite[Lemma 2.10]{kr-20}), that
\[
\big| (\Phi_2,q)_K\big| \le c\, \| u\|_{E_\beta^2(K)} \ \| q\|_{V_{1-\gamma}^2(K)}\le c'\, \| u\|_{E_\beta^2(K)} \ \| \phi\|_{V_{1-\gamma}^0(K)}\, .
\]
Analogously,
\[
\Big| \int_{\Gamma_2}\, 2\varepsilon_u(u)\, \cdot n\, \frac{\partial\bar{q}}{\partial n}\, dr\Big|
  \le c\, \| u\|_{E_\beta^2(K)} \ \| \phi\|_{V_{1-\gamma}^0(K)}
\]
and, consequently,
\[
\Big| \int p_2\, \bar{\phi}\, dx\Big| \le c\, \| u\|_{E_\beta^2(K)} \ \| \phi\|_{V_{1-\gamma}^0(K)}
\]
for all $\phi \in V_{1-\beta}^0(K)\cap V_{1-\gamma}^0(K)$. This proves (\ref{1al3}). From well-known regularity results
for elliptic equations in angles and cones it follows that $p_1 \in V_{\gamma}^1(K) \subset V_{\gamma-1}^0(K)$ if
$f\in V_\gamma^0(K)$, $g\in V_\gamma^1(K)\cap (V_{-\gamma}^1(K,\Gamma_2))^*$ and $h\in V_\gamma^{1/2}(\Gamma_2)$.
Thus $p=p_1+p_2 \in V_{\gamma-1}^0(K)$. The proof is complete. $\Box$

\subsubsection{Normal solvability}

First, we prove an estimate for the solution of the problem (\ref{par1}), (\ref{par2}) in the space $E_\beta^2(K,\Gamma_1)
\times {\cal V}_\beta^1(K)$ which is analogous to the estimates in \cite[Lemma 2.9]{kr-20} for the Dirichlet problem
and in \cite[Lemma 4.4]{kr-23} for the Neumann problem.

\begin{Le} \label{bl1}
Suppose that $u\in W^2({\cal D})$, $p\in W^1({\cal D})$ for every bounded subdomain ${\cal D}\subset K$, $\overline{\cal D}
\subset \overline{K}\backslash \{ 0\}$, and that $(u,p)$ is a solution of the problem {\em (\ref{par1}), (\ref{par2})}, where $s\not=0$.
If $u\in E_{\beta-1}^1(K,\Gamma_1)$, $p\in V_{\beta-1}^0(K)$, $x\cdot u \in (V_{2-\beta}^1(K,\Gamma_2))^*$,
$f\in V_\beta^0(K)$, $g\in V_\beta^1(K)\cap (V_{-\beta}^1(K,\Gamma_2))^*$ and $h\in E_\beta^{1/2}(\Gamma_2)$,
then $u\in E_\beta^2(K,\Gamma_1)$, $p \in V_\beta^1(K)$, $p|_{\Gamma_2}\in E_\beta^{1/2}(\Gamma_2)$ and
\begin{eqnarray} \label{1bl1}
&& \hspace{-2em}  \|  u\|_{V_{\beta}^2(K)} + |s|\, \| u\|_{V_{\beta}^0(K)}
  + \| p\|_{V_{\beta}^1(K)} + |s|^{1/4}\,  \| p\|_{V_{\beta}^0(\Gamma_2)} \nonumber \\ \nonumber
&& \le c\, \Big( \| f\|_{V_{\beta}^0(K)}  + \| g\|_{V_{\beta}^1(K)} + |s|\, \| g\|_{(V_{-\beta}^1(K,\Gamma_2))^*}
  + \| h\|_{V_{\beta}^{1/2}(\Gamma_2)}  \\
&& \quad + \ |s|^{1/4}\, \|  h\|_{V_{\beta}^0(\Gamma_2)} + \| u \|_{V_{\beta-1}^1(K)} + |s|^{1/2} \,  \| u\|_{V_{\beta-1}^0(K)}
  + \| p\|_{V_{\beta-1}^0(K)}\Big)
\end{eqnarray}
with a constant $c$ independent of $u,p,s$.
\end{Le}

Proof. Let $\chi_1$, $\chi_2$ be two times continuously differentiable functions in $K$ depending only on the variable
$\varphi=\mbox{arg}\, (x_1+ix_2)$ such that $\chi_1=1$ near $\varphi=0$, $\chi_2=1$ near $\varphi=\alpha$ and $\chi_1+\chi_2=1$.
Then $\chi_j u \in E_{\beta-1}^1(K)$, $\chi_j p\in V_{\beta-1}^1(K)$ and $x\cdot (\chi_1 u) \in (V_{2-\beta}^1(K))^*$.
The vector function  $\chi_j(u,p)$ satisfies the equalities
\begin{eqnarray*}
s\chi_j u - \Delta(\chi_j u) - \nabla\nabla\cdot (\chi_j u) + \nabla(\chi_j p) = f^{(j)}, \quad -\nabla\cdot(\chi_j u) = g^{(j)}
\end{eqnarray*}
in $K$, where $g^{(j)}= \chi_j g - u\cdot \nabla\chi_j$ and $f^{(j)}$ has the form
$f^{(j)} = \chi_j f + (L_1 u) \nabla\chi_j + (L_2\chi_j)\, u + p\nabla\chi_j$ with a homogeneous linear first order operator $L_1$
and a homogeneous linear second order operator $L_2$. Since $\partial_x^\alpha \chi_j(x) \le c\, r^{-|\alpha|}$ for $|\alpha|\le 2$, we obtain
\[
\| f^{(j)}\|_{V_\beta^0(K)} + \| g^{(j)}\|_{V_\beta^1(K)} \le c\, \Big( \| f\|_{V_\beta^0(K)}+ \| g\|_{V_\beta^1(K)} + \| u\|_{V_{\beta-1}^1(K)}
  + \| p\|_{V_{\beta-1}^0(K)}\Big)
\]
with a constant $c$ independent of $f,g,u$ and $p$. We estimate the norm of $g^{(j)}$ in $(V_{-\beta}^1(K,\Gamma_2))^*$. Let
$q\in V_{-\beta}^1(K,\Gamma_2)$. Then
\[
\int_K g^{(j)}\, q\, dx = \int_K \chi_j\, gq \, dx - \int_K q\, u\cdot\nabla \chi_j\, dx,
\]
where
\[
\Big| \int_K \chi_j g\, q\, dx \Big| \le \| g\|_{(V_{-\beta}^1(K,\Gamma_2))^*}\ \| \chi_j q\|_{V_{-\beta}^1(K)}
\]
and
\[
s \int_K q\, u\cdot\nabla \chi_j\, dx = \int_K q\, (f+\Delta u + \nabla\nabla\cdot u- \nabla p)\cdot\nabla \chi_j\, dx.
\]
Obviously,
\[
\Big| \int_K q\, f\cdot\nabla \chi_j\, dx\Big| \le \| f\cdot \nabla\chi_j\|_{V_{\beta+1}^0(K)}\, \| q\|_{V_{-\beta-1}^0(K)}
  \le c\, \| f\|_{V_\beta^0(K)}\, \| q\|_{V_{-\beta}^1(K)}\, .
\]
Since $q\nabla\chi_j=0$ on $\Gamma_1\cup \Gamma_2$, we obtain
\begin{eqnarray*}
\Big| \int_K q\nabla p \cdot \nabla \chi_j\, dx\Big| & = & \Big| \int_K p\, \nabla\cdot (q\nabla\chi_j)\, dx \Big|
  \le c\, \| p\|_{V_{\beta-1}^0(K)} \ \| q\|_{V_{-\beta}^1(K)} \, .
\end{eqnarray*}
Analogously,
\[
\Big| \int_K q\, (\Delta u + \nabla\nabla\cdot u)\cdot \nabla\chi_j\, dx\Big| \le c\, \| u\|_{V_{\beta-1}^1(K)} \ \| q\|_{V_{-\beta}^1(K)} \, .
\]
Thus,
\[
\Big| \int_K g^{(j)}\, q\, dx\Big| \le c\, \Big( \| g\|_{(V_{-\beta}^1(K,\Gamma_2))^*} + \| f\|_{V_\beta^0(K)} + \| u\|_{V_{\beta-1}^1(K)}
  + \| p\|_{V_{\beta-1}^0(K)}\Big) \, \| q\|_{V_{-\beta}^1(K)}
\]
for all $q\in V_{-\beta}^1(K,\Gamma_2)$. The vector function $\chi_1(u,p)$ satisfies the Dirichlet condition $\chi_1 u=0$ both on $\Gamma_1$
and $\Gamma_2$, whereas $\chi_2(u,p)$ satisfies the Neumann conditions
$-\chi_2pn + 2\varepsilon_n(\chi_2 u)=0$ on $\Gamma_1$ and $-\chi_2pn + 2\varepsilon_n(\chi_2 u)=h$ on $\Gamma_2$.
Using the estimates in \cite[Lemma 2.9]{kr-20} and \cite[Lemma 4.4]{kr-23}, we conclude that
\[
\| \chi_j u\|_{V_{\beta}^2(K)}+ |s|\, \| \chi_j u\|_{V_{\beta}^0(K)}+\| \chi_j p\|_{V_{\beta}^1(K)}+|s|^{1/4}\, \| \chi_j p\|_{V_{\beta}^0(\Gamma_2)}
\]
is majorized by the right-hand side of (\ref{1bl1}) for $j=1$ and $j=2$. Hence (\ref{1bl1}) is valid under the assumptions of the lemma. $\Box$ \\

Obviously, $x\cdot u \in V_{\beta-1}^0(K) \subset (V_{2-\beta}^1(K,\Gamma_2))^*$ for every $u\in E_\beta^2(K)$. Hence (\ref{1bl1}) is valid for
every solution $(u,p) \in E_\beta^2(K) \times {\cal V}_\beta^1(K)$ of the problem (\ref{par1}), (\ref{par2}).

\begin{Le} \label{bl2}
Suppose that $s\not=0$, the line $\mbox{\em Re}\, \lambda = \beta-1$ is free of solutions of the equation {\em (\ref{ev2})} and $\beta \not=
(k+\frac 12)\frac\pi\alpha$ for every integer $k$. Then every solution $(u,p)\in E_\beta^2(K,\Gamma_1)\times {\cal V}_\beta^1(K)$
of the problem {\em (\ref{par1}), (\ref{par2})} satisfies the estimate
\begin{eqnarray} \label{1bl2}
\| u\|_{E_\beta^2(K)} + \| p\|_{{\cal V}_\beta^1(K)} & \le & c\, \Big( \| f\|_{V_{\beta}^0(K)}  + \| g\|_{V_{\beta}^1(K)}
  + \| g\|_{(V_{-\beta}^1(K,\Gamma_2))^*} \nonumber\\
&& + \ \| h\|_{E_{\beta}^{1/2}(\Gamma_2)}
 + \| u \|_{W^1({\cal D})} +   \| p\|_{L_2({\cal D})}\Big),
\end{eqnarray}
where ${\cal D}$ is the set $\{ x\in K:\ \varepsilon < |x|<M \}$ with certain positive numbers $\varepsilon$ and $M$. The constant $c$
depends on $|s|,\varepsilon$ and $M$ but not on $u$ and $p$.
\end{Le}

Proof. As in the proof of \cite[Theorem 4.5]{kr-23}, we consider the cases that 1. $u$ and $p$ are zero outside a small neighborhood of the origin,
2. $u(x)$ and $p(x)$ are nonzero only for for large $|x|$ and 3. $u$ and $p$ are zero outside ${\cal D}$.

Case 1: If  $u(x)$ and $p(x)$ are zero for $|x|>2\varepsilon$, then it follows from Theorem \ref{at1} that
\[
\| u\|_{V_\beta^2(K)} + \| p\|_{V_\beta^1(K)} \le  c\, \Big( \| f-su\|_{V_{\beta}^0(K)}  + \| g\|_{V_{\beta}^1(K)}
  + \| h\|_{V_{\beta}^{1/2}(\Gamma_2)}\Big)
\]
Since
\[
\| u\|_{V_\beta^0(K)}\le \varepsilon^2 \, \| u\|_{V_\beta^2(K)} \quad\mbox{and}\quad \| p\|_{V_\beta^0(\Gamma_2)}
  \le c\, \varepsilon^{1/2} \, \| p\|_{V_\beta^1(K)}\, ,
\]
we obtain (\ref{1bl1}) if $\varepsilon$ is sufficiently small.

Case 2: If $u(x)$ and $p(x)$ are zero for $|x|<M/2$, then
\begin{equation} \label{2bl2}
\| u\|_{E_{\beta-1}^1(K)} \le c\, M^{-1}\, \| u\|_{E_\beta^2(K)} \, .
\end{equation}
Let $\gamma$ be a real number, $\beta < \gamma \le \beta+\frac 12$, such that the interval $\beta \le \lambda\le \gamma$
does not contain numbers of the form $(k+\frac 12)\frac\pi\alpha$ with integer $k$.
Then it follows from Lemma \ref{al3} that $p=p_1+p_2$, where
$p_1,p_2 \in V_\beta^1(K)$, $p_1$ satisfies (\ref{estp1}) and $p_2$ satisfies the estimate (\ref{1al3}).
Since $p(x)$ is equal to zero for $|x|<M$, we get
\begin{eqnarray} \label{3bl2} \nonumber
&& \| p\|_{V_{\beta-1}^0(K)} \le \| p_1\|_{V_{\beta-1}^0(K)} + M^{\beta-\gamma}\, \| p_2\|_{V_{\gamma-1}^0(K)}
  \le   c\, \Big( \| f\|_{V_\beta^0(K)} + \| g\|_{V_\beta^1(K)} \\
&& \qquad + \ \| g\|_{(V_{-\beta}^1(K,\Gamma_2))^*} + \| h\|_{V_\beta^{1/2}(\partial K)}
  + M^{\beta-\gamma}\, \| u \|_{E_\beta^2(K)}\Big).
\end{eqnarray}
If $M$ is sufficiently large, then the estimates (\ref{1bl1}), (\ref{2bl2}) and (\ref{3bl2}) imply (\ref{1bl2}).

Case 3: If $u$ and $p$ are zero outside ${\cal D}$, then the estimate (\ref{1bl2}) follows directly from (\ref{1bl1}).

Combining the estimates for the solution in the cases 1--3 as in the proof of \cite[Theorem 2.1]{kr-20}, we obtain the estimate
(\ref{1bl2}) for arbitrary solutions $(u,p)\in E_\beta^2(K,\Gamma_1)\times {\cal V}_\beta^1(K)$ of the problem
(\ref{par1}), (\ref{par2}). $\Box$ \\

Since every bounded subset of $E_\beta^2(K,\Gamma_1)\times {\cal V}_\beta^1(K)$ is precompact in $W^1({\cal D}) \times L_2({\cal D})$,
the following assertion holds as a consequence of Lemma \ref{bl2}.

\begin{Co} \label{bc1}
Suppose that $s\not=0$, the line $\mbox{\em Re}\, \lambda = \beta-1$ is free of solutions of the equation {\em (\ref{ev2})} and $\beta \not=
(k+\frac 12)\frac\pi\alpha$ for every integer $k$. Then the range of the operator $A_{\beta,s}$ is closed and its kernel has finite dimension.
\end{Co}

\subsection{A regularity assertion for strong solutions}

In order to apply Lemma \ref{bl1}, we need the following assertion on the function $x\cdot u$.

\begin{Le} \label{bl3}
Suppose that $s\not=0$ and that $(u,p) \in E_\beta^2(K,\Gamma_1) \times {\cal V}_\beta^1(K)$ is a solution of the problem
{\em (\ref{par1}), (\ref{par2})}. Assume furthermore that $f\in V_\gamma^0(K)$, $g\in V_\gamma^1(K)\cap
(V_{-\gamma}^1(K,\Gamma_2))^*$ and $h\in V_\gamma^{1/2}(\Gamma_2)$, where $\beta\le \gamma\le \beta+\frac 12$, and that
all numbers $(k+\frac 12)\frac\pi\alpha$ with integer $k$ lie outside the interval $\beta \le \lambda \le \gamma$. Then
$x\cdot u \in (V_{2-\gamma}^1(K,\Gamma_2))^*$.
\end{Le}

Proof. By (\ref{par1}), we have $sx\cdot u = x\cdot (f+ \Delta u + \nabla\nabla\cdot u - \nabla p)$. Let $q\in V_{2-\gamma}^1(K,\Gamma_2)$. Then
\[
\Big| \int_K q\, x\cdot f\, dx\Big| \le \| f\|_{V_\gamma^0(K)}\ \| q\|_{V_{1-\gamma}^0(K)}
  \le \|f\|_{V_\gamma^0(K)}\ \| q\|_{V_{2-\gamma}^1(K)}\, .
\]
By Lemma \ref{al3}, we have $p\in V_{\gamma-1}^0(K)$. Since $x\cdot n=0$ on $\partial K\backslash \{ 0\}$, we get
\begin{eqnarray*}
&& \Big|\int_K q\, x\cdot \nabla(\nabla\cdot u -p)\, dx\Big| = \Big| \int_K (g+p)\, \nabla\cdot (qx)\, dx\Big| \\
&&  \le \| g+p\|_{V_{\gamma-1}^0(K)} \ \| xq\|_{V_{1-\gamma}^1(K)} \le (\| g\|_{V_\gamma^1(K)}+ \| p\|_{V_{\gamma-1}^0(K)}) \
  \| q\|_{V_{2-\gamma}^1(K)} \, .
\end{eqnarray*}
Finally,
\[
\int_K q\, x\cdot \Delta u \, dx = -\sum_{j=1}^2 \int_K \nabla u_j \cdot \nabla(x_j q)\, dx
  + \int_{\Gamma_1} q\, x \cdot \frac{\partial u}{\partial n}\, dr.
\]
Here,
\[
\Big| \int_K \nabla u_j \cdot \nabla(x_j q)\, dx\Big| \le  \| u_j\|_{V_{\gamma-1}^1(K)}) \   \| q\|_{V_{2-\gamma}^1(K)}
  \le  \| u\|_{E_\beta^2(K)}) \   \| q\|_{V_{2-\gamma}^1(K)}
\]
and
\[
\Big| \int_{\Gamma_1} q\, x \cdot \frac{\partial u}{\partial n}\, dr\Big| \le \Big\| r^{\gamma-1/2} \frac{\partial u}{\partial n}\Big\|_{L_2(\Gamma_1)} \
  \| r^{-\gamma+3/2} q\|_{L_2(\Gamma_1)} \le \Big\| \frac{\partial u}{\partial n}\Big\|_{E_\beta^{1/2}(\Gamma_1)}\   \| q\|_{V_{2-\gamma}^1(K)}
\]
since $E_\beta^2(K) \subset V_{\gamma-1}^0(K)$ and $E_\beta^{1/2}(\Gamma_1) \subset V_{\gamma-1/2}^0(\Gamma_1)$ if $\beta\le \gamma\le \beta+\frac 12$
(see Subsection 3.1). This proves the lemma. $\Box$ \\

Now the following lemma can be proved analogously to \cite[Lemma 2.9]{kr-16} and \cite[Lemma 4.7]{kr-23}.

\begin{Le} \label{bl4}
Suppose that $(u,p) \in E_\beta^2(K,\Gamma_1) \times {\cal V}_\beta^1(K)$ is a solution of the problem {\em (\ref{par1}), (\ref{par2})},
where $s\not=0$ and
\begin{equation} \nonumber \label{1bl3}
f \in V_\beta^0(K)\cap V_\gamma^0(K), \ g \in V_\beta^1(K) \cap V_\gamma^1(K)\cap (V_{-\beta}^1(K,\Gamma_2))^* \cap (V_{-\gamma}^1(K,\Gamma_2))^*, \
\end{equation}
and $h\in E_\beta^{1/2}(\Gamma_2)\cap E_\gamma^{1/2}(\Gamma_2)$. Furthermore, let one of the following two conditions be satisfied.
\begin{itemize}
\item[(i)] $\beta<\gamma$ and the interval $\beta \le \lambda \le \gamma$
  does not contain numbers of the form $(k+\frac 12)\frac \pi\alpha$ with integer $k$,
\item[(ii)] $\gamma<\beta$ and the strip $1-\beta \le \mbox{\em Re}\, \lambda\le 1-\gamma$ contains no solutions of the equation {\em (\ref{ev2})}.
\end{itemize}
Then $u\in E_\gamma^2(K,\Gamma_1)$ and $p\in {\cal V}_\gamma^1(K)$.
\end{Le}

Proof. Let $\zeta$ be a two times continuously differentiable function on $\overline{K}$
with compact support which is equal to one in a neighborhood of the origin, and let $\eta=1-\zeta$.

(i) Suppose that $\beta<\gamma$. Then $\zeta u \in  E_\gamma^2(K,\Gamma_1)$ and $\zeta p \in  {\cal V}_\gamma^1(K,\Gamma_2)$.
We assume in addition that $\gamma\le \beta+\frac 12$. Then $\eta u \in E_\beta^2(K) \subset E_{\gamma-1}^1(K)$,
$\eta p \in V_{\gamma-1}^0(K)$ (by Lemma \ref{al3}) and $\eta x\cdot u \in (V_{2-\gamma}^1(K,\Gamma_2))^*$ (by Lemma \ref{bl3}).
Obviously,
\[
(s-\Delta)(\eta u) -\nabla\nabla\cdot(\eta u) +\nabla(\eta p) = \eta f + L_\eta u + p\nabla \eta \in V_\gamma^0(K)
\]
since $L_\eta u = \eta(\Delta u + \nabla\nabla\cdot u) -\Delta(\eta u)-\nabla\nabla\cdot(\eta u)$ is a sum of terms 
$a_{\alpha,\beta}\partial_x^\alpha \eta \, \partial_x^\beta u$, where $|\alpha|+|\beta|=2$ and $|\beta|\le 1$.
Furthermore,
\[
 -\nabla\cdot (\eta u) = \eta g-u\cdot\nabla\eta \in V_\gamma^1(K)\cap (V_{-\gamma}^1(K,\Gamma_2))^*,
\]
$\eta u=0$ on $\Gamma_1$ and $-\eta pn + 2\varepsilon_N(\eta u) \in E_\gamma^{1/2}(\Gamma_2)$. Using Lemma \ref{bl1},
we obtain $\eta (u,p) \in E_\gamma^2(K,\Gamma_1) \times {\cal V}_\gamma^1(K)$. This proves that
$(u,p) \in E_\gamma^2(K,\Gamma_1) \times {\cal V}_\gamma^1(K)$ if $\beta\le \gamma\le \beta+\frac 12$. By induction we obtain
the same result for $\beta+\frac k2 \le \gamma \le \beta+\frac{k+1}2$, $k=1,2,\ldots$.

(ii) Suppose that $\beta>\gamma$. Then $\eta u \in  E_\gamma^2(K,\Gamma_1)$ and $\eta p \in  {\cal V}_\gamma^1(K)$.
If in addition $\gamma \ge \beta-2$, then $\zeta u \in E_\beta^2(K)\subset V_\gamma^0(K)$. Thus, it can be easily shown
by means of Theorem \ref{at1} that $\zeta (u,p) \in V_\beta^2(K)\times V_\beta^1(K)$. This proves that
$(u,p) \in E_\gamma^2(K,\Gamma_1) \times {\cal V}_\gamma^1(K)$ if $\beta-2\le \gamma\le \beta$. By induction,
the same result holds for $\gamma <\beta-2$. $\Box$ \\

\subsubsection{Existence of strong solutions}

Let again $\lambda_1$ be the smallest positive solution of the equation (\ref{ev2}).
We consider the operator $A_{\beta,s}$ defined by (\ref{operator}).

\begin{Le} \label{bl5}
The operator $A_{\beta,s}$ is injective if $-\frac{\pi}{2\alpha} < \beta < 1+\mbox{\em Re}\, \lambda_1$, $\mbox{\em Re}\, s\ge 0$ and $s\not= 0$.
\end{Le}

Proof. Suppose that $(u,p) \in E_\beta^2(K,\Gamma_1) \times {\cal V}_\beta^1(K)$ is a solution of the problem (\ref{par1}), (\ref{par2})
with the data $f=0$, $g=0$ and $h=0$.

1) If $0\le \beta\le 1$, then $(u,p)\in E_0^1(K,\Gamma_1)\times (L_2(K)+\stackrel{\circ}{V}\!{}_0^1(K,\Gamma_2))$.
Furthermore, $(u,p)$ satisfies the identity (\ref{1al2}) with $F=0$ for all $v\in E_0^1(K,\Gamma_1)$ and the equation
$\nabla \cdot u=0$ in $K$ (see Lemma \ref{al2}). Using Theorem \ref{at1}, we obtain $u=0$ and $p=0$.

2) Suppose that $-\frac{\pi}{2\alpha} < \beta <0$. Then Lemma \ref{bl4} implies $u \in E_0^2(K,\Gamma_1)$ and
$p\in {\cal V}_0^1(K)$. From 1) it follows that $(u,p)=(0,0)$.

3) Suppose $1<\beta< 1+\mbox{Re}\, \lambda_1$. Then Lemma \ref{bl4} implies $u \in E_1^2(K,\Gamma_1)$ and
$p\in {\cal V}_1^1(K)$. From 1) it follows that $(u,p)=(0,0)$. $\Box$

\begin{Th} \label{bt1}
Suppose that $s\not=0$, $\mbox{\em Re}\, s \ge 0$, $\alpha<\pi$ and
\begin{equation} \label{Ialpha}
\max\Big(-\frac{\pi}{2\alpha},1-\lambda_1\Big) < \beta < \min\Big(\frac{\pi}{2\alpha},1+\lambda_1\Big)
\end{equation}
Then the operator $A_{\beta,s}$ is bijective, and every solution $(u,p) \in E_\beta^2(K,\Gamma_1) \times {\cal V}_\beta^1(K)$
of the problem {\em (\ref{par1}), (\ref{par2})} satisfies the estimate {\em (\ref{1bt1})}
with a constant $c$ independent of $f,g,h$ and $s$.
\end{Th}

Proof. First note that the $\beta$-interval (\ref{Ialpha}) is not empty since $\lambda_1>1/2$ for $\alpha<\pi$.
From our assumptions it follows that the line $\mbox{Re}\, \lambda = \beta-1$ is free of solutions of the equation (\ref{ev2}) and $\beta \not=
(k+\frac 12)\frac\pi\alpha$ for every integer $k$. By Corollary \ref{bc1} and Lemma \ref{bl5}, the range of the operator $A_{\beta,s}$ is closed
and its kernel is trivial. We prove the existence of solutions $(u,p) \in E_\beta^2(K,\Gamma_1) \times {\cal V}_\beta^1(K)$
for arbitrary $f\in V_\beta^0(K)$, $g\in V_\beta^1(K)\cap (V_{-\beta}^1(K,\Gamma_2))^*$, $h\in E_\beta^{1/2}(\Gamma_2)$.
For this, we consider the following two cases separately.

1) $\beta \in I_\alpha \cap [0,1]$, where $I_\alpha$ denotes the interval (\ref{Ialpha}).
By Corollary \ref{ac3}, there exists a unique weak solution $(u,p)\in E_0^1(K,\Gamma_1)\times (L_2(K)+V_0^1(K,\Gamma_2))$
of the problem (\ref{par1}), (\ref{par2}). Let $\zeta$ be the same cut-off function as in the prof of Lemma \ref{bl4} and
$\eta=1-\zeta$. Then $\zeta(u,p) \in V_0^1(K) \times L_2(K)$ and
\[
-\Delta(\zeta u) -\nabla\nabla\cdot(\zeta u)+ \nabla(\zeta p) = \zeta(f-su) -L_\eta u +p\nabla\zeta \in V_\gamma^0(K),
\]
where $L_\eta$ is the same differential operator as in the proof of Lemma \ref{bl4}. Furthermore
\[
-\nabla\cdot(\zeta u)=\zeta g -u\cdot\nabla\zeta \in V_\gamma^1(K), \quad -\zeta p n+ 2\varepsilon_n(\zeta u)\in V_\gamma^{1/2}(\Gamma_2)
\]
for arbitrary $\gamma\ge \beta$. From Theorem \ref{at1} it follows that $\zeta(u,p) \in V_1^2(K) \times V_1^1(K)$
Since the strip  $0\le \mbox{Re}\, \lambda \le 1-\beta$ is free of solutions of the equation (\ref{ev2}), we conclude from Theorem \ref{at1}
that $\zeta u \in V_\beta^2(K)$ and $\zeta p \in V_\beta^1(K)$.
Furthermore,
\[
(s-\Delta)(\eta u) -\nabla\nabla\cdot(\eta u)+ \nabla(\eta p) = \eta f +L_\eta u +p\nabla\eta \in V_\gamma^0(K),
\]
$-\nabla\cdot(\eta u) \in V_\gamma^1(K)\cap (V_{-\gamma}^1(K,\Gamma_2))^*$ and $-\eta p n+ 2\varepsilon_n(\eta u)\in E_\gamma^{1/2}(\Gamma_2)$
for arbitrary $\gamma\ge \beta$. Since moreover $x\cdot u \in L_2(K) \subset (V_1^1(K))^*$, we have
$\eta x\cdot u \in (V_{2-\gamma}^1(K,\Gamma_2))^*$ for $\gamma\le \beta\le 1$. Obviously, $\eta u\in E_{-1}^1(K,\Gamma_1)$
and $\eta p \in V_{-1}^0(K)$. Thus, Lemma \ref{bl1} implies $\eta(u,p) \in E_0^2(K,\Gamma_1)\times {\cal V}_0^1(K)$.
Using Lemma \ref{bl4}, we conclude that $\eta(u,p) \in E_\beta^2(K,\Gamma_1)\times {\cal V}_\beta^1(K)$
and, consequently, $(u,p) \in E_\beta^2(K,\Gamma_1)\times {\cal V}_\beta^1(K)$.

2) Suppose that $\alpha<\frac \pi 2$ (i.e., $\lambda_1>1$) and $\beta\in I_\alpha \backslash [0,1]$. In addition to
the above assumption on $f,g,h$, we assume that $f\in V_\gamma^0(K)$, $g\in V_\gamma^1(K)\cap (V_{-\gamma}^1(K,\Gamma_2))^*$
and $h\in E_\gamma^{1/2}(\Gamma_2)$ with a certain $\gamma \in I_\alpha\cap [0,1]$. As was just shown, there exists
a solution $(u,p) \in E_\gamma^2(K,\Gamma_1) \times {\cal V}_\gamma^1(K)$ of the problem (\ref{par1}), (\ref{par2}).
Using  Lemma \ref{bl4}, we conclude that $(u,p) \in E_\beta^2(K,\Gamma_1) \times {\cal V}_\beta^1(K)$.
Since the range of the operator $A_{2,\beta}$ is closed, it follows that the problem (\ref{par1}), (\ref{par2})
is solvable in $E_\beta^2(K,\Gamma_1) \times {\cal V}_\beta^1(K)$ for arbitrary
$f\in V_\beta^0(K)$, $g\in V_\beta^1(K)\cap (V_{-\beta}^1(K,\Gamma_2))^*$, $h\in E_\beta^{1/2}(\Gamma_2)$.
This proves the bijectivity of the operator $A_{\beta,s}$.

In particular, it follows that the solution $(u,p)$ satisfies (\ref{1bt1}) if $\mbox{Re}\, s\ge 0$ and $|s|=1$. Here the constant
$c$ is independent of $\mbox{arg}\, s$. Suppose that $\mbox{Re}\, s\ge 0$, $s\not=0$ and $(u,p) \in E_\beta^2(K)\times {\cal V}_\beta^1(K)$
is a solution of the problem (\ref{par1}), (\ref{par2}). We define $v(x)=u(|s|^{-1/2}x)$, $q(x)=|s|^{-1/2}p(|s|^{-1/2}x)$,
$F(x)=|s|^{-1}f(|s|^{-1/2}x)$, $G(x)=|s|^{-1/2}g(|s|^{-1/2}x)$ and \\ $H(x)= |s|^{-1/2}h(|s|^{-1/2}x)$. Then $(v,q)$ is a solution
of the equation
\[
A_{\beta,|s|^{-1}s}\, (v,q) = (F,G,H).
\]
If $\beta$ satisfies the inequalities (\ref{Ialpha}), then
\[
\| v\|_{E_\beta^2(K)} + \| p\|_{V_\beta^1(K)} \le c\, \Big( \| F\|_{V_\beta^0(K)} + \| G\|_{V_\beta^1(K)} +
 \| G\|_{V_{-\beta}^1(K,\Gamma_2))^*} + \| H\|_{E_\beta^{1/2}(\Gamma_2)}\Big)
\]
with a constant $c$ indepe4ndent of $s$. This implies (\ref{1bt1}). $\Box$

\section{The time-dependent problem in \boldmath $K$}

Now, we consider the problem (\ref{stokes1a})--(\ref{stokes3a}).
As in Sections 2--4, we assume that $K$ is a two-dimensional angle with opening $\alpha$.

\subsection{Weighted Sobolev spaces in $K\times{\Bbb R}_+$}

Let $Q=K\times {\Bbb R}_+ = K\times (0,\infty)$ and $\Gamma = ({\partial K}\backslash \{ 0\})\times {\Bbb R}_+$.
We denote by $W_\beta^{2,l}(Q)$ the weighted Sobolev space of all functions ${\mathfrak u}(x,t)$ on $Q$ with finite norm
\[
\| {\mathfrak u}\|_{W_\beta^{2,1}(Q)} = \Big( \int_0^\infty  \big( \| {\mathfrak u}(\cdot,t) \|^2_{V_\beta^{2}(K)}
  + \| \partial_t {\mathfrak u}(\cdot,t) \|^2_{V_\beta^0(K)}\big)\, dt \Big)^{1/2}.
\]
The space $\stackrel{\circ}{W}\!{}_\beta^{2,1}(Q)$ is the subspace of all ${\mathfrak u}\in W_\beta^{2,1}(Q)$ satisfying
the condition ${\mathfrak u}(x,0)$ for $x\in K$.  We also need the trace space $W_\beta^{1/2,1/4}(\Gamma_2)$.
The norm in this space is defined as
\begin{eqnarray*}
\| {\mathfrak h}\|_{W_\beta^{1/2,1/4}(\Gamma_2)} & = & \Big( \int_0^\infty \| {\mathfrak h}(\cdot,t)\|^2_{V_\beta^{1/2}(\Gamma_2)}\, dt  \\
&& \quad +\ \int_\Gamma r^{2\beta} \int_0^\infty \int_0^t \frac{|{\mathfrak h}(x,t)-{\mathfrak h}(x,t-\tau)|^2}{\tau^{3/2}}\,
   d\tau\, dt\, dx\Big)^{1/2}.
\end{eqnarray*}
We consider the Laplace transforms. Let $H_\beta^2(K)$ be the space of holomorphic functions
$u(x,s)$ for $\mbox{Re}\, s>0$ with values in $E_\beta^2(K)$ for which the norm
\[
\| u\|_{H_\beta^2(K)} = \sup_{\gamma>0} \Big( \frac 1i \int_{\mbox{\scriptsize Re}\, s=\gamma}  \big( \| u(\cdot,s)\|^2_{V_\beta^2(K)}
  + |s|^2  \| u(\cdot,s)\|^2_{V_\beta^0(K)}\big)\, ds  \Big)^{1/2}
\]
is finite.  Analogously, $H_\beta^{1/2}(\Gamma_2)$ is the space of holomorphic functions
$h(x,s)$ for $\mbox{Re}\, s>0$ with values in $E_\beta^{1/2}(\Gamma_2)$ for which the norm
\[
\| h\|_{H_\beta^{1/2}(\Gamma_2)} = \sup_{\gamma>0} \Big( \frac 1i \int_{\mbox{\scriptsize Re}\, s=\gamma}   \big( \| h(\cdot,s)
  \|^2_{V_\beta^{1/2}(\Gamma_2)} + |s|^{1/2}  \| h(\cdot,s)\|^2_{V_\beta^0(\Gamma_2)}\big)\, ds  \Big)^{1/2}
\]
is finite. The following lemma was stated in \cite[Proposition 3.4]{kozlov-89}. The proof is essentially the same as for nonweighted
spaces in \cite[Theorem 8.1]{av-64}.

\begin{Le} \label{cl1}
The Laplace transform realizes isomorphisms between the spaces $\stackrel{\circ}{W}\!{}_\beta^{2,1}(Q)$ and
$W_\beta^{1/2,1/4}(\Gamma)$ on one side and the spaces $H_\beta^2(K)$ and $H_\beta^{1/2}(\partial K)$ on the other side.
\end{Le}

\subsection{Existence and uniqueness of solutions}

As a consequence of Theorem \ref{bt1}, we obtain the following assertion.

\begin{Th} \label{ct1}
Suppose that $\alpha<\pi$, ${\mathfrak f}\in L_2\big( {\Bbb R}_+,V_\beta^0(K)\big)$, ${\mathfrak g}\in L_2\big({\Bbb R}_+,V_\beta^1(K)\big)$,
$\partial_t g\in L_2\big({\Bbb R}_+,(V_{-\beta}^1(K,\Gamma_2))^*)\big)$, ${\mathfrak h} \in W_\beta^{1/2,1/4}(\Gamma_2)$ and that ${\mathfrak g}(x,0)=0$
for $x\in K$. Furthermore, we assume that $\beta$ satisfies the inequalities {\em (\ref{Ialpha})}. Then there exists a unique
solution $({\mathfrak u},{\mathfrak p}) \in \stackrel{\circ}{W}\!{}_\beta^{2,1}(Q) \times L_2\big({\Bbb R}_+, V_\beta^1(K)\big)$
of the problem {\em (\ref{stokes1a})--(\ref{stokes3a})} satisfying the estimate
\begin{eqnarray*}
&& \hspace{-1em}\| {\mathfrak u}\|_{W_\beta^{2,1}(Q)} + \|{\mathfrak p}\|_{L_2({\Bbb R}_+,V_\beta^1(K))} +
  \| {\mathfrak p}\|_{W_\beta^{1/2,1/4}(\Gamma_2)}  \\
&&  \hspace{-1em}\le c\, \Big( \| {\mathfrak f}\|_{L_2({\Bbb R}_+,V_\beta^0(K))} + \| {\mathfrak g} \|_{L_2({\Bbb R}_+,V_\beta^1(K))}
  +  \| {\mathfrak g}_t\|_{L_2({\Bbb R}_+,(V_{-\beta}^1(K,\Gamma_2))^*)} + \| {\mathfrak h}\|_{W_\beta^{1/2,1/4}(\Gamma_2)}\Big)
\end{eqnarray*}
with a constant $c$ independent of ${\mathfrak f}$, ${\mathfrak g}$ and ${\mathfrak h}$.
\end{Th}

P r o o f.
Let $f$ and $g$ and $h$ be the Laplace transforms (with respect to the variable $t$) of ${\mathfrak f}$, ${\mathfrak g}$ and ${\mathfrak h}$,
respectively. For arbitrary $s\not=0$, $\mbox{Re}\, s\ge 0$, there exists a uniquely determined solution
$(u,p)(\cdot,s) \in E_\beta^2(K)\times {\cal V}_\beta^1(K)$ of the problem (\ref{par1})
satisfying the estimate
\begin{eqnarray*}
&& \| u(\cdot,s)\|^2_{V_\beta^2(K)} + |s|^2\, \| u(\cdot,s)\|^2_{V_\beta^0(K)}
  + \| p(\cdot,s)\|^2_{V_\beta^1(K)} + |s|^{1/2} \, \| p(\cdot,s)\|_{V_\beta^0(\Gamma_2)} \\
&& \le c\, \Big( \| f(\cdot,s)\|^2_{V_\beta^0(K)} + \| g(\cdot,s) \|^2_{V_\beta^1(K)} + |s|^2\, \| g(\cdot,s)\|^2_{(V_{-\beta}^1(K,\Gamma_2)^*)} \\
&& \qquad   + \ \| h(\cdot,s)\|^2_{V_\beta^{1/2}(\Gamma_2} + |s|^{1/2}\, \| h(\cdot,s)\|_{V_\beta^0(\Gamma_2)}\Big)
\end{eqnarray*}
with a constant $c$ independent of $s$. Integrating over the line $\mbox{Re}\, s = \gamma$ and taking the supremum
with respect to $\gamma>0$, we get the assertion of the theorem. \hfill $\Box$

\end{document}